\documentclass[12pt]{article} 

\usepackage{amsmath,amsxtra,amssymb,latexsym, amscd}
\usepackage[mathscr]{eucal}

\begin{document}

\title{ \huge\bf Gr\"obner Basis\\Convex Polytopes and \\
Planar Graph
 }

\def\1{\rule{0cm}{0cm}} \def\qd{\rule{3mm}{3mm}} \def\BB{$\bullet$}
\renewcommand{\arraystretch}{1.25}
\renewcommand{\theequation}{\thesectn.\arabic{equation}}
\def\sce{\setcounter{equation}{0}}  \newcounter{sectn} 
\newcounter{sbsect}
\def\sect#1{\addtocounter{section}{1}\sce\setcounter{sbsect}{0}%
        \renewcommand{\thesectn}{\thesection}\1\smallskip\\
        {\1\hspace{-2em}\large\bf\thesectn.\qquad #1\smallskip\par}}
\def\subsect#1{\addtocounter{sbsect}{1}\sce%
        
\renewcommand{\thesectn}{\thesection:\Alph{sbsect}}\1\smallskip\\
        {\bf\1\hspace{-1.5em}\thesectn.\qquad #1\smallskip\par}}
\newtheorem{Theorem}{THEOREM} \newtheorem{Lemma}[Theorem]{LEMMA}
\newtheorem{Corollary}[Theorem]{COROLLARY}
\def\thm#1#2{\be{Theorem}{\lb{#1} #2}} \def\LEM#1#2{\BE{Lemma}{\LB{#1} 
#2}}
\def\COR#1#2{\BE{Corollary}{\LB{#1} #2}}
\def\proof{\bigskip\noindent {\sc Proof:}\qquad}
\def\REM{\1\smallskip\par\noindent{\bf REMARK:}\qquad }
\def\qed{\hfill$\quad$\qd\medskip\\} \def\ds{\displaystyle}
\def\LB#1{\label{#1}} \def\BE#1#2{\begin{#1} #2 \end{#1}}
\def\EQ#1#2{\BE{equation}{\LB{#1} #2}} \def\ARR#1#2{\BE{array}{{#1} 
#2}}
\def\DES#1{\BE{description}{#1}} \def\QT#1{\BE{quote}{#1}}
\def\ENUM#1{\BE{enumerate}{#1}} \def\ITM#1{\BE{itemize}{#1}}
 \def\COM#1{\par\noindent{\bf COMMENT:\quad\sl #1}\par\noindent}
%\def\COM#1{}   % so comments don't print!!
%%%%%%%%%%%%%%%%%%%%%%%%%%%%%%%%%%%%%%%%%%%%%%%%%%%%%%%%%%%%%%%%%%%%%%
\def\mapsfrom{\hbox{$\;{\leftarrow}\kern-.15em{\mapstochar}\:\:$}}
\def\vv{\kern.344em{\rule[.18ex]{.075em}{1.32ex}}\kern-.344em}
\def\RE{\mbox{\rm I\kern-.21em R}} \def\CX{\mbox{\rm \vv C}}
\def\imp{\Rightarrow} \def\emb{\hookrightarrow} 
\def\wk{\rightharpoonup}
\def\rd{\dot{\1}} \def\d{\cdot} \def\+{\oplus} \def\x{\times}
\def\<{\langle} \def\>{\rangle} \def\o{\circ} \def\at#1{\Bigr|_{#1}}
\def\cd{\partial} \def\grad{\nabla} \def\L{\left} \def\R{\right}
\def\bx{\mathbf{x}} \def\by{\mathbf{y}} \def\bS{\mathbf{S}}
\def\I{{\cal I}} \def\A{{\cal A}} \def\D{{\cal D}}\def\bc{\mathbf{c}}
\def\bx{\mathbf{x}} \def\by{\mathbf{y}} \def\bS{\mathbf{S}}
\def\H{{\mathcal H}} \def\U{{\cal U}} \def\D{{\cal 
D}}\def\bc{\mathbf{c}}
%w%%%%%%%%%%%%%%%%%%%%%%%%%%%%%%%%%%%%%%%%%%%%%%%%%%%%%%%%%%%%%%%%%
\def\eq{equation} \def\de{differential \eq} \def\pde{partial \de}
\def\sol{solution} \def\pb{problem} \def\bdy{boundary} 
\def\fn{function}
\def\dde{delay \de} \def\ev{eigenvalue}
\def\R{\mathbb R}
\def\C{\mathbb C}
\author{
{\bf Dang Vu Giang}\\
Hanoi Institute of Mathematics\\
18 Hoang Quoc Viet, 10307 Hanoi, Vietnam\\
{\footnotesize          e-mail: $\<$dangvugiang@yahoo.com$\>$}\\
\1\\
}
\maketitle 
{\footnotesize
\noindent {\bf Abstract.}    Using the Gr\"obner basis of an ideal generated by a family of polynomials we prove that every planar graph is 4-colorable. Here we also use the fact that the complete graph of 5 vertices is not included in any planar graph.

\medskip


\medskip

}

\bigskip

\par\noindent Consider a planar graph $G$ with $n$ vertices. Between two any vertices there is at most one edge. We want to color the vertices of G by 4 colors such that no two vertices connected by an edge are colored the same way. We let $\zeta =i\in \mathbb{C}$ be the 4th primitive root of unity $\left( {{i}^{4}}=1 \right).$ Let ${{x}_{1}},{{x}_{2}},\cdots ,{{x}_{n}}$ be variables representing the distinct vertices of the graph $G.$ Each vertex is to be assigned one of  colors $1,i,{{i}^{2}},{{i}^{3}}.$ This means that $x_{k}^{4}-1=0$ for $k=1,2,\cdots ,n.$ Moreover, if ${{x}_{j}}$ and ${{x}_{k}}$ are connected then $x_{j}^{3}+x_{j}^{2}{{x}_{k}}+{{x}_{j}}x_{k}^{2}+x_{k}^{3}=0.$  Let $I$ be the ideal generated by ${{P}_{k}}=x_{k}^{4}-1$ for $k=1,2,\cdots ,n$ and ${{Q}_{j,k}}=x_{j}^{3}+x_{j}^{2}{{x}_{k}}+{{x}_{j}}x_{k}^{2}+x_{k}^{3}$ if ${{x}_{j}}$ and ${{x}_{k}}$ are connected. We should prove that $I\ne \mathbb{C}\left[ {{x}_{1}},{{x}_{2}},\cdots ,{{x}_{n}} \right].$ In fact, the complete graph of 5 vertices is not included in any planar graph. Not every ${{Q}_{j,k}}$ is a generator of  I. Therefore, the set of generators of I is not so large. Now we induct on $n$. Assume that every planar graph of $n-1$ vertices is $4-$colorable. Now consider a planar graph G of $n$ vertices. If some vertex of G has degree less than 4 then G is clearly, $4-$colorable. Hence, we assume that vertex of G has degree larger than 3. Let ${{x}_{n}}$ be the vertex of smallest degree and let ${{x}_{1}},{{x}_{2}},\cdots ,{{x}_{k}}$ be vertices connected with ${{x}_{n}}.$ Here, $k$ is the degree of the vertex ${{x}_{n}}.$ Thus,  $k\ge 4.$ Without loss of generality we assume that ${{x}_{1}},{{x}_{2}},\cdots ,{{x}_{k}}$ are colored by 4 distinct colors. Now assume for the sake of a contradiction that ${{I}_{n}}=\mathbb{C}\left[ {{x}_{1}},{{x}_{2}},\cdots ,{{x}_{n}} \right].$  
That is a real ideal  ${{I}_{n-1}}$ of $\mathbb{C}\left[ {{x}_{1}},{{x}_{2}},\cdots ,{{x}_{n-1}} \right]$ together  with ${{P}_{n}}=x_{n}^{4}-1$ and ${{Q}_{j,n}}=x_{n}^{3}+x_{n}^{2}{{x}_{j}}+{{x}_{n}}x_{j}^{2}+x_{j}^{3}$   for $j=1,2,\cdots ,k$ generates the whole ring $\mathbb{C}\left[ {{x}_{1}},{{x}_{2}},\cdots ,{{x}_{n}} \right].$ We can write $1=\left( x_{n}^{4}-1 \right){{\varphi }_{0}}+\sum\limits_{j=1}^{k}{\left( x_{n}^{3}+x_{n}^{2}{{x}_{j}}+{{x}_{n}}x_{j}^{2}+x_{j}^{3} \right){{\varphi }_{j}}+f\left( {{x}_{1}},{{x}_{2}},\cdots ,{{x}_{n-1}} \right)}$ where $f\in {{I}_{n-1}}$ and ${{\varphi }_{j}}\in \mathbb{C}\left[ {{x}_{1}},{{x}_{2}},\cdots ,{{x}_{n}} \right]$  for $j=0,1,2,\cdots ,k.$

{\footnotesize  
\bigskip

\par\noindent{\bf Acknowledgement.} Deepest appreciation is extended towards the NAFOSTED  (the National Foundation for Science and Techology Development in Vietnam) for the financial support.}

\end{document}